\documentclass[12pt]{amsart}
\usepackage[utf8]{inputenc}
\usepackage{%
  mathtools,  
  amssymb,    
  amsthm,     
  amsfonts,   
  thmtools,   
  graphicx,
  float,
  color,
  xcolor,
  tikz,
  tikz-cd,
  mathrsfs,   
  etoolbox    
}

\definecolor{darkred}{rgb}{0.75,0,0}
\def\customcitecolor{darkred}
\def\customlinkcolor{darkred}

\usepackage[%
    colorlinks,
    citecolor=\customcitecolor,%
    linkcolor=\customlinkcolor,%
    urlcolor=\customlinkcolor%
]{hyperref}

\usepackage[capitalise,nameinlink,noabbrev]{cleveref}

\usepackage[margin=1in]{geometry}

\usepackage[final,nopatch=footnote]{microtype}

\theoremstyle{definition}
\newtheorem{theorem}{Theorem}[section]
\numberwithin{theorem}{section} 
\numberwithin{equation}{section} 

\usepackage{mfirstuc}
\newcommand{\capitalizename}[1]{\makefirstuc{#1}}

\newcommand{\defthm}[1]{%
  \newtheorem{#1}[theorem]{\capitalizename{#1}}%
}

\newcommand{\defthms}[1]{%
  \forcsvlist{\defthm}{#1}%
}

\defthms{%
  answer,assumption,claim,conjecture,construction,corollary,
  counterexample,definition,digression,discussion,example, 
  examples,exercise,fact,goal,idea,intuition,lemma,
  motivation,notation,note,problem,proposition,question,remark,setup,
  slogan,strategy,terminology,upshot,warning%
}

\Crefname{problem}{Problem}{Problems}

\Crefname{question}{Question}{Questions}

\newcommand{\deftextcommand}[1]{%
  \expandafter\providecommand\csname #1\endcsname{\mathrm{#1}}%
}
\newcommand{\deftextcommands}[1]{%
    \forcsvlist{\deftextcommand}{#1}%
}

\deftextcommands{ab,alg,an,ann,Aut,BG,BGL,Bl,BO,BP,BSL,BSO,BSp,BSU,BU,can,cd,cdh,cl,coBar,codim,codom,coeq,coev,cof,cofib,coker,colim,coim,cone,conj,const,coTor,cyc,diag,Desc,dg,Disc,disc,dual,eff,EKL,End,eq,ess,et,Et,EU,ev,Ex,ex,Exc,Ext,fib,Fix,Fl,fppf,fpqc,Frac,Frob,Fun,Gal,gen,GL,gp,Gr,gr,GW,Her,Ho,hocofib,hocolim,hofib,holim,Hom,id,Idem,im,incl,Ind,ind,inj,Inn,Inv,inv,iso,Jac,KGL,kgl,KH,KO,ko,KQ,kq,KR,KSp,KU,ku,Lan,Map,map,MGL,MO,Mor,mor,MSL,MSO,MSp,MSU,MU,mult,MUP,Nm,ob,obj,op,Orb,ord,Out,perf,Perm,PGL,Pin,pr,pre,Proj,proj,prom,PSL,quot,Ran,rank,Res,RO,sep,sgn,SH,sig,Sing,SL,SO,soc,Sp,Span,Spec,Spin,spn,Sq,st,Stab,SU,supp,Supp,Syl,syl,Sym,syn,SYT,TC,td,Th,THH,Tor,Tot,TP,TR,Tr,tr,triv,univ,var,veff,vir,vol,Wel,Wr}

\deftextcommands{Ab,Aff,Alg,Ani,Bimod,CAlg,Cat,CDGA,CG,CGWH,Ch,CMon,coAlg,Coh,CommRing,ConjSub,coMod,Cor,Corr,CoSh,CRing,CW,Field,Fin,FinSet,Gpd,Grp,Grpd,Grph,Kan,Kar,LMod,Mfld,Mod,NAlg,Ouv,Perf,Poset,Pr,Pre,PSh,PShv,qCat,QCoh,Rep,Ring,RMod,sAb,Set,SH,Sh,Shv,Sm,Sp,Spc,Spectra,sPre,sSet,sShv,Stack,Sub,Top,Tors,Var,Vect}

\newcommand{\defblackboardletter}[1]{%
  \expandafter\providecommand\csname #1\endcsname{\mathbb{#1}}
}
\newcommand{\defblackboardletters}[1]{%
  \forcsvlist{\defblackboardletter}{#1}%
}

\defblackboardletters{A,C,F,P,Q,R,Z}

\providecommand{\RP}{{\mathbb{R}\text{P}}}

\RequirePackage{bbm}

\RequirePackage{pict2e,picture}
\makeatletter
\DeclareRobustCommand{\DDelta}{{\mathpalette\bb@Delta\relax}}
\newcommand{\bb@Delta}[2]{%
  \begingroup
  \sbox\z@{$\m@th#1\Delta$}%
  \dimendef\Dht=6 \dimendef\Dwd=8
  \setlength{\Dwd}{\wd\z@}%
  \setlength{\Dht}{\ht\z@}%
  \begin{picture}(\Dwd,\Dht)
  \put(0,0){$\m@th#1\Delta$}
  \put(.42\Dwd,.7\Dht){\line(10,-26){.25\Dht}}
  \end{picture}%
  \endgroup
}

\usepackage{ytableau}

\newcommand{\MgnbarXb}
{\overline{M}_{g,n}(X,\beta)}
\newcommand{\MzeronbarXb}
{\overline{M}_{0,n}(X,\beta)}

\title[The Evolution of Enumerative Geometry]{The Evolution of Enumerative Geometry:  A Narrative from Classical Problems to Enriched Invariants} 

\author{Candace Bethea}
\author{Thomas Brazelton}
\date{\today}

\begin{document}

\begin{abstract}

Enumerative geometry, the art and science of counting geometric objects satisfying geometric conditions, has seen a resurgence of activity in recent years due to an influx of new techniques that allow for \emph{enriched} computations. This paper offers a historical survey of enumerative geometry, starting with its classical origins and real counterparts, to new advances in  quadratic enrichment. We include a brief survey of the paradigm shift initiated by Gromov-Witten theory, whose impact can be seen in recent results in quadratically enriched enumerative geometry. Finally, we conclude with a brief overview of emerging directions including random and equivariant enumerative geometry. 
\end{abstract}

\maketitle

\section{Introduction}

Enumerative geometry is, colloquially speaking, the study of counting geometric things. When a mathematician asks ``how many $X$ satisfy property $Y$'' they are often met with the answers zero or infinitely many. Enumerative geometry deals with questions whose answers lie between these, where one has a finite amount of data to latch onto and study. Famous examples in this area include the Apollonian problem of finding the eight circles tangent to three on the plane, Salmon and Cayley's computation that there are 27 lines on a smooth cubic surface, and Kontsevich's recursive formula for the finite numbers $N_d$ of degree $d$ rational curves interpolating $3d-1$ points on the projective plane.

The power of enumerative geometry comes from the statement of \emph{conservation of number}, which says that the answer to the problem is well-defined and independent of initial parameters, provided the parameters are generically chosen. For example any two smooth cubic surfaces have 27 lines, or very pedantically, the set of lines on any smooth cubic surface has cardinality 27. When we ask a specific enumerative geometry problem, e.g., finding the flexes on a smooth planar cubic, we are naturally handed a \emph{finite set} of solutions. Conservation of number says that given two instances of the problem, the finite sets of solutions are the same in number, i.e., when we count solutions with multiplicity, we obtain the same value. To that end we posit the following:

\textbf{Pseudo-definition:} We say a theory of enumerative geometry is \emph{enriched} if its solutions are sets equipped with additional structure that generalize multiplicity (we typically refer to this generalization as \emph{weight} or \emph{mass}), and it has a conservation of number statement that the weighted sums of these structures agree for any two generic instantiations of the problem.

In classical enumerative geometry, solutions are weighted by their multiplicity. In this note, we explore a few instances of enriched enumerative geometry, including:
\begin{itemize}
    \item \emph{Signed real enumerative geometry}: We work over the real numbers, and each solution is weighted not only by its multiplicity but also by a sign $+$ or $-$, where a negative sign captures a local change in orientation. The sum of these signed contributions is conserved so long as the sum is finite.
    \item \emph{Gromov--Witten theory}: We generally work over the complex numbers, though considerable work has been done to study Gromov--Witten invariants over other fields. Each solution is weighted by a finite group of automorphisms, and in general the sum of the reciprocals of the orders of automorphism groups of solutions is conserved.
    \item \emph{Quadratically enriched enumerative geometry}: We work over an arbitrary field $k$, and each solution is weighted by a (virtual) symmetric bilinear form over $k$. The solution is conserved in the group completion of the monoid of isomorphism classes of symmetric bilinear forms, the \emph{Grothendieck--Witt ring} $\GW(k)$.
    \item \emph{Random enumerative geometry}: We work over the real numbers, and instantiate a problem according to some notion of randomness. The finite set of solutions is instead replaced by an \emph{expected value} of the number of solutions.\footnote{This doesn't quite fit the pseudo-definition of enrichment above, but it fits nicely into the story we're telling here.}
    \item \emph{Equivariant enumerative geometry}: We are again over the complex numbers, working over an object with some finite symmetry group $G$ acting on the moduli space of potential solutions. Our solution sets are naturally $G$-sets, and conservation says that two sets of solutions are in $G$-equivariant bijection.
\end{itemize}

The purpose of these conference proceedings is to expose readers to the development of some of these ideas and how they flow together. A tremendous amount of  detail, nearly all, is missing from these notes, so it is \emph{not intended to be a technical introduction or a thorough historical introduction to any one of these collections of work}. As J.F. Adams once wrote, this work should be viewed as ``an essay in machine appreciation; it is not intended to qualify the reader for a mechanic's certificate'' \cite{adams_infinite_1978}. In particular, this article should be supplemented by the many cited references for anyone interested in legitimately learning the mathematics discussed here. To that end, we have omitted certain definitions and results and skipped certain parts of historical narrative, instead focusing  on specific things according to our own interests and within the context of this conference proceedings. Nevertheless we hope this might serve as interesting reading to anyone who has ever heard the phrase ``enriched enumerative geometry'' and wondered what it was all about.

\textbf{Acknowledgements}: The authors would like to thank Park City Math Institute and the organizers of the summer research program in Motivic Homotopy Theory that took place during July 2025. We would like to thank the following people for helpful comments and correspondence that tremendously improved earlier drafts of this paper: Jim Bryan, Andrés Jaramillo Puentes, Steve Kleiman, César Lozano Huerta, Jake Levinson, Stephen McKean, and Nicolas Michel. Candace Bethea was supported by National Science Foundation award DMS-2402099. Thomas Brazelton was supported by National Science Foundation award DMS-230324.

\section{Classical enumerative geometry}

While enumerative geometry is often said to have begun in antiquity with the work of the great Greek geometers, perhaps its two most important components (algebraic closure and conservation of number) did not crystallize until the 19th century. Throughout the Renaissance, original Greek texts saw direct Latin translations, rather than translation by way of the Arabic translations from the Islamic Golden Age centuries prior. The wave of interest in ancient Greek culture at this time was felt throughout many disciplines and art forms, and mathematics was not immune: a renewed focus was placed upon studying classical Greek geometry in the style of the ancients. Moreover the now mythic status of these geometers of antiquity led mathematicians to devote a great deal of effort to work out exactly what mathematics may have been contained in the many texts lost to time. A key example of this is Viète's 1600 text \emph{Apollonius Gallus} \cite{viete_apollonius_1600} in which he attempts to reconstruct the arguments from the lost books in Apollonius of Perga's foundational treatment on conics. While mathematicians both during the Renaissance and today are familiar with Euclid's \emph{Elements}, many believe that the loftiest peaks of Greek geometry are lost to time. Throughout the mathematical literature of antiquity, \emph{Elements} is contrasted with Euclid's lost work \emph{Porisms}, in which ``higher geometry'' was said to be developed. Mathematicians including Maclaurin and Newton have developed theories about what may have been contained in this text, and the tradition of postulating about the Porisms continues into the modern day.

We begin our story of classical enumerative geometry in France at the beginning of the 19th century. At the École Polytechnique, geometry was taught alongside analysis to civil engineers, politicians, military officers, as well as future mathematicians. Analysis, in contrast to algebra or geometry, enjoyed a rare position of prestige, viewed even outside mathematics as a tremendous educational endeavor critical to the education of anyone desiring to rise the ranks of Napoleonic France. Gaspard Monge, a professor of geometry at École Polytechnique, published his treatise on descriptive geometry\footnote{Descriptive geometry should be distinguished from \emph{analytic geometry} in that the latter makes reference to coordinates, while the former attempts to be as coordinate-free as possible.} in 1799 --- this is demarcated by Coolidge as the beginning of the resurrection of projective geometry from its Greek origins. 
Monge was an incredibly renowned geometer in his own right, but his influence is perhaps best felt through the generation of students he graduated, and the subsequent geometric revolution that their ideas inspired. He held his position until he was ousted by Laplace (who detested Monge and devalued his mathematics by extension) in 1815.

We focus on two key figures who studied at École Polytechnique during this period of time:  Jean-Victor Poncelet and Michel Chasles. Both shared a dream of establishing a \emph{pure geometry}, which would bring geometry onto the same footing as algebra and analysis. Pure geometry would be a language, akin to algebra or analysis, in which computations could be carried out. The so-called \emph{ideal solutions} would exist in direct analogy to how imaginary numbers exist as solutions of real equations \cite{chemla_fragments_2024} and would allow for general geometric statements which treat many cases at the same time. Both Poncelet, and later Chasles, viewed Monge's descriptive geometry as a stepping stone towards this goal, although they had disagreements in its scope. 

One of Poncelet's core principles which he wished to include in pure geometry was his notion of \emph{continuity} -- this held that certain geometric invariants\footnote{The application to enumerative geometry which would come later is the behavior of \emph{projective} invariants under continuous changes, although Poncelet studies many metric invariants in his book.}
were invariant under small parametrized changes. As a motivating (if not mildly ahistorical) example, consider the following:

\begin{problem}[The Circles of Apollonius]
\label{problem:apollonius}
Given three circles on the plane, how many are tangent to all three?\footnote{This problem is attributed to Apollonius of Perga (who Geminus of Rhodes nicknamed the \emph{excellent geometer}), 
although no original writing on this problem exists. The attribution comes from nearly half a millennium later, in work of Pappus of Alexandria, which was resurrected and popularized by Vi\`{e}te \cite{viete_apollonius_1600}.}
\end{problem}
The answer to this problem (provided we interpret circles and planes in a projective way over the complex numbers, which the reader may correctly object is ahistorical) is eight.
In this context, we may interpret Poncelet's principle of continuity as indicating the following: given three initial circles and the eight tangent circles, by slightly perturbing the equation of any of our initial circles, the tangent circles grow or shrink or translate in order to accommodate tangency, but no new solutions are created and no existing solutions destroyed. We invite the reader to think about this as homotopy invariance, and indeed this can be reinterpreted as the statement that a general algebraic section of $\mathcal{O}_{\mathbb{P}^3}(2)^{\oplus 3}$ is transverse to the zero section at eight points,\footnote{%
While $\mathbb{P}^5$ is a moduli space of conics, enumerative questions dealing with circles leverage projective three-space as a moduli space of circles. This can be described of the linear subspace $\mathbb{P}^3 \subset \mathbb{P}^5$ of those conics passing through the two points $[1:i:0]$ and $[1:-i:0]$ (called the \emph{circle points} by Chasles). The locus of circles tangent to a fixed one forms a quadric hypersurface in $\mathbb{P}^3$, and by Bézout's theorem, three generic such hypersurfaces meet transversely in eight distinct points.
}
and this fact is invariant under a small perturbation of the section.\footnote{The invariance of this statement boils down to the observation that locally, transversality is the statement that a particular matrix has nonzero determinant. This property is preserved under any very small deformation of the entries in the matrix.}

Michel Chasles' studies at École Polytechnique were briefly interrupted by his conscription into Napoleon's forces, after which he tried his hand in finance before returning to mathematics. In 1837 he published his \emph{Aperçu Historique}, considered one of the great historical texts on geometry, but also one of the most questionable. While Poncelet frequently had to teach analysis and viewed pure geometry as an independently valuable field of geometry, he considered coordinate-dependent arguments as a distraction from the pure intuition which must underly geometry; Poncelet was not an undivisive figure in his endeavors, getting into public arguments with figures like Cauchy and Gergonne. Chasles, on the other hand, drifted into flawed historiography in his efforts to champion pure geometry. His attempts to recreate Euclid's \emph{Porisms} cast them as part of a direct mathematical genealogy which passed through Monge's descriptive geometry and arrived at his own work \cite{Michel_Smadja_2022}. Like Poncelet, he situated pure geometry in direct contrast to analytic geometry, but pursued its development in a very different way.

Chasles' theory of \emph{characteristics} provided brilliant insight, and set the stage for what would become Schubert calculus. As an example, consider the following:

\begin{problem}
    Given $i$ points and $5-i$ lines on the plane, how many conics pass through these points and are tangent to these lines?
\end{problem}
Working with $\mathbb{P}^5$ as a space of conics, Bézout's theorem tells us to compute the degree of the hypersurfaces corresponding to passing through a point or being tangent to a line, which are one and two, respectively. The answer then appears to be $2^{5-i}$, a statement implied by de Jonquière's prior work \cite{de_jonquieres_sur_1866}. 
In hindsight, we can appeal to pole-polar duality to see that this cannot be correct, as it implies in the dual projective plane that there are 32 conics passing through five points. The problem with this is a defect in $\mathbb{P}^5$ as a moduli space for conics -- it contains a locus of doubled lines and line pairs which can contribute degenerate solutions to results that attempt to leverage Bézout's theorem on this space.\footnote{Explicitly, there is a 2-dimensional family of conics described by the square of an equation of a line. Geometrically, these are lines considered with multiplicity two - in particular any point of intersection between a doubled line and another curve is (technically speaking) a point of tangency, and contributes to a count of conics satisfying a given tangency condition. One has to throw out this contribution in order to obtain the desired honest geometric count.}

Chasles' insight was that, given a family of conics, the condition $\mu$ of passing through a point and $\nu$ of being tangent to a line should both be considered, and any condition on conics is a linear combination of these. The number of conics passing through $i$ points and tangent to $5-i$ lines can then be computed as $2^{\text{min}\{i,5-i\}}$. Using his theory of characteristics, Chasles counted the smooth conics tangent to five given conics as 3264, correcting an erroneous computation of Steiner which fell victim to the same issue creating the number 32 above. In contemporary language, Chasles' computations can be thought of as occurring in the Chow ring of the moduli space of complete conics, where $\mu$ and $\nu$ are the pullback of the hyperplane class to the projective plane and its dual, respectively \cite[\S8.2]{3264}. We refer the reader to \cite{Kleiman-Chasles} for a more detailed account of this story.

In the next generation of geometers, a student of Chasles named Hieronymous Zeuthen built upon the foundation of pure geometry that had been initiated. Unlike Poncelet and Chasles who preferred a pure geometry in which pictures did not guide intuition, Zeuthen visualized singularities and intersection multiplicity in a dynamic way, and carried out explicit computations with algebraic equations in order to study geometry. Zeuthen's survey paper on enumerative geometry for the German encyclopedia of mathematics, later translated into French together with Pieri, highlights the usage of Poncelet's conservation of number in enumerative geometry, and is still considered one of the great texts in enumerative geometry.

Hermann Schubert was also deeply influenced by Chasles' perspectives on geometry, publishing his thesis on the three-dimensional analog of the Circles of Apollonius problem. Inspired by Chasles and by Poncelet's principle of continuity, Schubert advocated for what he called the \emph{principle of conservation of number}. This held that the numerical answer to an enumerative problem was independent of the initial parameters chosen (e.g., there are three circles tangent to \emph{any} three on the plane, or 27 lines on \emph{any} cubic surface), provided they are chosen suitably generically. While Poncelet's idea of continuity communicated that relationships between geometric figures were invariant under small parametrized changes, Schubert's perspective was decidedly more static --- he viewed the principle of conservation of number as a computational technique: an enumerative theorem may be proven by picking parameters and then solving this instantiation of the given problem.

With this principle in hand, Schubert began developing a calculus to enumerate geometric objects in three-dimensions satisfying certain conditions (these and many of their higher-dimensional generalizations are now called \emph{Schubert problems}). These ranged from relatively simple problems like counting the two lines meeting four lines in three-space, to counting the number of quadric hypersurfaces tangent to nine generic quadric hypersurfaces in three-space. Schubert called his mathematics the \emph{calculus of conditions} (we now call it \emph{Schubert calculus}). Kleiman's overview of the developments of Hilbert's 15th problem and Schubert calculus give a significantly more detailed historical overview as well as an exposition of beautiful examples of Schubert calculus in action \cite{kleiman-hilbert-15}. 

Despite the fact that all his computations have proven over time to be correct, Schubert provided no semblance of a proof that his calculus of conditions was valid. This drew the ire of mathematicians including Zeuthen, Halphen, Study, and Kohn. 
The shaky bedrock upon which this theory was built led to Hilbert including \emph{Schubert calculus} as the fifteenth of his eponymous problems at the turn of the century.

\subsection{Resolving Hilbert's 15th problem}

Schubert's calculus of conditions, Chasles's theory of characteristics, and much of 19th century enumerative geometry, lacked a rigorous foundation. This framework emerged from work of van der Waerden in 1929 \cite{van_der_waerden_topologische_1930}, leveraging simplicial cohomology, a tool coming from algebraic topology, and work of Lefschetz \cite{browder_problem_1976}. The idea was the following: rather than study some algebraic locus (say, a hypersurface representing some geometric condition) we should instead investigate its associated cohomology class. The condition of a circle being tangent to a fixed circle, for instance, is the class $2h\in H^2(\mathbb{CP}^3;\Z)$, and changing the fixed circle does not affect the cohomology class. To carry out a computation which combines conditions, we want to compute the number of points in the intersection, which is obtained by multiplying these classes using the ring structure on singular cohomology. To compute the number eight in the Circles of Apollonius problem for instance, we compute that the number of circles tangent to three of them is $2^3 h^3 \in H^6(\mathbb{CP}^3;\Z)$, which after integrating is equal to $8$. One may say that the Chow ring is the home for Schubert's computations, as it is the home for intersection theory of algebraic varieties, however most of the spaces over which these classical computations are carried out were cellular varieties, hence the cycle class map from the Chow ring to singular cohomology is an isomorphism, and the computation is agnostic as to whether we view it as occurring in Chow groups or cohomology.\footnote{Viewing computations as Chow-valued can mislead us in attempts to enrich enumerative computations to other settings --- for instance attempting to build equivariant enumerative geometry using equivariant Chow groups produces answers valued in the complex representation ring of the group, which is not sufficient to recover the $G$-action on the solutions themselves (see for instance \cite[Example~5.22]{EEG}).} Viewing these computations as occurring in singular cohomology is not only historically sound, it is pedagogically useful here -- the observation that an enumerative computation can be carried out via a characteristic class is crucial in efforts to enrich classical enumerative geometry to other contexts. It also provides an elegant resolution to the principle of continuity as discussed, that a characteristic class computed by a section depends only on the homotopy class of the section -- i.e. it is invariant under \textit{continuous} deformations, not just algebraic ones.

\subsection{Parameterized problems: towards thinking about stacks} Following the development of theories like singular cohomology, the difference between counting solutions to an enumerative problem and computing a homology class became more or less indistinguishable. Indeed many enumerative problems whose rigorous solutions eluded 19th century geometers can be now thought of as easy computations in the cohomology ring of a certain moduli space. Schubert calculus has a very different meaning today than it did a century ago; instead of referring to the procedure of solving and enumerating linear intersection problems it now means (to many) the study of cohomology of homogeneous spaces under a simple Lie group.


Many moduli spaces one wishes to study, however, cannot be accurately captured as an algebraic variety or manifold. Often this is due to the moduli space attempting to describe objects which are \emph{parametrized} in some sense. A motivating problem is the following classical interpolation problem:

\begin{problem}
\label{problem:count_plane_curves}
How many rational degree $d$ curves pass through $3d-1$ points in general position on $\mathbb{P}^2$?
\end{problem}

The development of \emph{stacks} in the 1960's introduced a mathematical foundation that refined the existing study of moduli spaces for these sorts of questions. This led to a resurgence of enumerative geometry as well as a shift of focus towards Gromov-Witten theory and related theories, which we discuss more in \Cref{section:Gromov-Witten}.

\section{Real counts in enumerative geometry}\label{section:real}

The study of algebraic geometry over the reals is much more complicated than over the complex numbers due to the lack of algebraic closure. In enumerative algebraic geometry this is felt in the failure of conservation of number. Consider for instance the classical count of 27 lines on a smooth complex cubic surface \cite{cayley_triple_1849}.

\begin{theorem}[{\cite{schlafli_attempt_1858}}]
A real smooth cubic surface can have 3, 7, 15, or 27 real lines.
\end{theorem}

One perspective to take on this is that the moduli of smooth complex cubic surfaces is connected, and the incidence variety of cubic surfaces equipped with a line ramifies only over the singular cubic surfaces, hence any smooth cubic surface contains the same number of lines over $\mathbb{C}$. Over the reals, the moduli of smooth real cubic surfaces is disconnected, having five connected components corresponding to the number of real lines and real tritangents \cite[\S23]{segre_non-singular_1942}.

Indeed this type of behavior is the primary complicating factor in studying smooth objects in real algebraic geometry --- over the complex numbers a codimension one discriminant locus has real codimension two, whereas over the reals it has real codimension one, chopping up a space into many connected components. As an example to keep in mind when contemplating complexity, every smooth planar curve of degree $d$ over $\mathbb{C}$ looks identical to a topologist, whereas over the real numbers the number of connected components in the moduli space of real smooth planar curves of degree $d$ grows exponentially in $d^2$ \cite{orevkov_asymptotic_2003}; the classification of their possible shapes is part of Hilbert's 16th problem, and is completely open for $d\ge8$.

Perhaps one of the earliest cases of real counts to enumerative problems is one we have already seen -- this is the Circles of Apollonius (\Cref{problem:apollonius}). The problem of counting circles tangent to three over the reals was well-studied by Vi\`ete at the turn of the 17th century in his famous book \emph{Apollonius Gallus} \cite{viete_apollonius_1600}, 
after which it remained a consistent topic of study, enjoying solutions by Newton, Gergonne, and many others.\footnote{For an overview of solutions to this problem see \cite{court_problem_1961}.}
For the purposes of this narrative, however, we view this topic as lying more firmly in classical Euclidean and projective geometry than in enumerative algebraic geometry.

In the late 1930's, Beniamino Segre (a former student of Corrado Segre) had been studying the geometry and topology of real algebraic varieties, particularly intersections, via their limiting behavior. After fleeing fascism in Italy, he was interned in England, and it was during this time he worked on perhaps his most well-known work, which was his treatise on non-singular cubic surfaces \cite{segre_non-singular_1942}. This book begins by studying the behavior of lines under a degeneration from a cubic surface to a union of three planes. This topological perspective, in contrast to the more algebraic one, allows for a more careful geometric analysis of cubic surfaces over the real numbers, as highlighted by Zariski in his review of Segre's book. In particular, given a real line on a cubic surface, and a hyperplane containing that line, the hyperplane cuts the cubic surface at a residual conic, intersecting the line at two points (Segre called these \textit{parabolic points}). Interchanging pairs of parabolic points gives an involution of any real line, and lines are called \textit{hyperbolic} or \textit{elliptic} corresponding to whether the fixed locus of the involution is real or complex, respectively \cite[\S27]{segre_non-singular_1942}. This leads to the following theorem.

\begin{theorem}
On a real smooth cubic surface, we have that
\begin{equation}\label{eqn:hyp-ellip}
\begin{aligned}
    \#\{\text{hyperbolic lines}\} - \#\{\text{elliptic lines}\} = 3.
\end{aligned}
\end{equation}
\end{theorem}

This can be proven in a number of ways, for instance Benedetti and Silhol proved that a real cubic surface inherits a $\Pin^-$ structure whose modulo four reduction can distinguish the two types of lines at the level of homology \cite{benedetti_rm_1995}. This can also be proven from the perspective of open Gromov-Witten theory, see for instance \cite{Solomon-thesis,horev_open_2012}. A perspective we discuss here is what one may call \emph{absolute Euler classes}.

\begin{definition}[{\cite{okonek_intrinsic_2014}}]
\label{def:absolute-euler-class}
Let $X$ be a closed real manifold of dimension $n$, and $V \to X$ a real topological vector bundle which is \emph{relatively oriented}\footnote{C.f.~\Cref{def:euler-class-GW}, note that for real line bundles being a square is the same as being trivial.} in the sense that the line bundle $\Hom(\det TX, \det V)$ admits a trivialization $\theta$. Then $V$ admits an Euler class $e(V,\theta)$ which depends on $\theta$ up to a sign. We call $|e(V,\theta)|$ the \emph{absolute Euler class} of $V$.
\end{definition}

Analogous to how a Chern class of a corank zero bundle is Poincar\'e dual to the vanishing locus of a generic section, the absolute Euler class provides a \emph{signed count} of the zeros of a generic section. For example the absolute Euler class of $\Sym^3\mathcal{S}^\ast \to \Gr_\R(2,4)$ is equal to $3$, and the local index at a line is equal to $+1$ or $-1$ corresponding to whether the line is hyperbolic or elliptic, respectively. From this perspective \Cref{eqn:hyp-ellip} can be reinterpreted as a formula arising from an absolute Euler class, and it becomes clear how to generalize this to other settings. In \cite{okonek_intrinsic_2014,finashin_abundance_2013} the authors compute the absolute Euler class for the analogous symmetric bundles over the Grassmannian of lines in higher-dimensional space, proving for instance a lower bound of $(2n-1)!!$ for the number of real lines on a general smooth real hypersurface of degree $2n-1$ in $\mathbb{P}^{n+1}$. A beautiful analog of this result for the 240 $(-1)$-curves on a real degree one del Pezzo can be found in \cite{finashin_two_2021}, a spiritually similar result for bitangents to real algebraic curves can be found in \cite{blomme_bitangents_2024}, and a signed count of rational curves on a generic real K3 surface can be found in \cite{MR3384445}.

\subsection{Signed real enumerative geometry}\label{section:real-signed-counts} The lesson to be learned from Segre's work is the following: while conservation of number may break over $\R$, a certain \emph{signed} count of solutions may remain invariant, and this sign can encode beautiful information about the local geometry. This is the jumping off point for quadratically enriched enumerative geometry, which we touch on more in \Cref{section:quadratically_enriched}. Before doing this, we discuss another key appearance of this idea, which is that of \emph{Welschinger invariants}.

Recall that the number $N_d$ of complex rational curves interpolating $3d-1$ points on $\mathbb{P}^2$ in \Cref{problem:count_plane_curves} is invariant of the position of the points, provided they are in general position. As we might expect, this same statement fails in the context of real curves -- for instance through 8 generic points in $\RP^2$ there can be 8, 10, or 12 real rational cubics interpolating them \cite{degtyarev_topological_2000}. Groundbreaking work of Welschinger tells us we should count the cubics interpolating these points weighted by a sign. More explicitly:

\begin{definition}
\label{def:welschinger-invariant} If $C$ is a rational curve of degree $d$, we define its \emph{Welschinger invariant} to be $\Wel(C) = (-1)^n$, where $n$ is the number of isolated points of $C$.\footnote{Recall a real point on a real rational curve is said to be \emph{isolated} if its directions of tangency form a complex conjugate pair.} This is also sometimes called the \emph{mass} of the curve.
\end{definition}

The remarkable theorem is the following:

\begin{theorem}[{\cite{Welschinger-GAFA}}]
\label{thm:welschinger} For $3d-1 = n_1 + 2n_2$ points on $\RP^2$, $n_1$ of them real and $n_2$ pairs of complex conjugate points, we have that the quantity

\begin{align*}
    W_{d,n_1} = \sum_{\substack{C\text{ real deg }d \\ \text{through these points}}} \Wel(C).
\end{align*}
is independent of the choice of points (provided they are chosen generically).
\end{theorem} 

\section{The Emergence of Gromov--Witten theory}\label{section:Gromov-Witten}

A survey of enumerative geometry would be remiss to omit a discussion of Gromov--Witten theory, currently one of the most active areas of enumerative geometry. The goal of this section is a departure from prior sections; this section primarily serves to introduce high-level aspects of the subject that are most relevant to current trends in quadratically enriched enumerative geometry, discussed in Section \ref{section:quadratically_enriched}. This section does not seek to serve as a comprehensive historical review or a complete technical introduction to the subject, which would be beyond the scope of this work. Considering this warning to the reader, we include  technical references introducing the subject in varying levels of detail \cite{13/2-notes, fulton-pandharipande, kock-vainsencher, simon-rose-notes}. 


Motivated by symplectic geometry and rigidity questions, Gromov introduced the study of pseudo-holomorphic curves and their existence in the presence of a symplectic structure on an almost complex smooth manifold \cite{gromov-pseudoholomorphic-curves}. In 1988, groundbreaking work of Witten introduced topological sigma models as another means of studying maps from the Riemann sphere to an almost complex manifold \cite{witten-sigma-models}. These beautiful new perspectives on the study of maps from spheres to complex manifolds arguably ignited the modern study of Gromov--Witten invariants in symplectic and algebraic geometry and topology, which we focus on in this section.

Our perspective on Gromov--Witten theory is the algebraic study of rational curves in algebraic varieties, a question which has appeared in various forms throughout this article: 

\begin{problem}\label{problem:rational_curve_count}
How many rational curves lie on a smooth complex projective variety? 
\end{problem}
When $X$ is $\mathbb{CP}^2$, a rational degree $d$ curve on $\mathbb{CP}^2$ is a map $f\colon \mathbb{CP}^1\to \mathbb{CP}^2$ given by $f([s:t]) = [f_0(s,t): f_1(s,t): f_2(s,t)]$ where the $f_i$ are degree $d$ homogeneous polynomials. This is the natural algebraic analogue of the question of studying maps from the Riemann sphere to a complex manifold, which is a simple to state and notoriously difficult to answer. 

The starting point for Problem \ref{problem:rational_curve_count} is for $\mathbb{CP}^2$, which is the question of how many degree $d$ rational plane curves exist for any $d$. With no additional restrictions, there are infinitely many for each $d$. We may ask the same question after imposing the condition that we seek to count rational curves that pass through a prescribed number of points in general position in $\mathbb{CP}^2$. A dimension argument shows that the space of degree $d$ rational, nodal plane curves is $3d-1$ dimensional, which leads us to computing the finite solutions $N_d$ to \Cref{problem:count_plane_curves}. Many well-known variations on  this problem are considered classical, and they have historically captured the attention of algebraic geometers independent of developments Gromov--Witten theory. A celebrated result of Caporaso-Harris leverages the geometry of the Severi varieties parameterizing nodal plane curves of degree $d$ and geometric genus $g$ to give a recursive formula enumerating the number of nodal plane curves of degree $d$ passing through an appropriate number of points \cite{caporaso-harris}.

One of the early successes of Gromov--Witten theory was the proof of a recursive formula in $d$ for the number of degree $d$ rational plane curves, denoted $N_d$, given by Kontsevich and Manin \cite[5.2.1]{kontsevich-manin}. In rapid succession, Gromov--Witten theory was developed symplectically and algebraically  \cite{kontsevich1, ruan-tian1, ruan-tian2, mcduff-salamon, behrend-manin, behrend-GW}, leading to a flurry of results in enumerative geometry that had previously been unattainable using existing methods. Though these results are too numerous to name individually, we mention a few \cite{equivariant-GW, vakil-rational-surfaces, yau-zaslow}. 

Given a smooth, projective, complex algebraic variety $X$  and $\beta\in H_2(X, \mathbb{Z})$, Gromov--Witten theory studies stable maps from nodal curves of arbitrary genus to $X$. The moduli stack of genus $g$, $n$-marked stable maps, denoted $\MgnbarXb$ 
is a proper Deligne-Mumford stack of finite type \cite[Section 1.3.1 p.3]{kontsevich1}. Note when $X = \mathbb{CP}^2$, $\beta = d\cdot \ell$ where $d\geq 1$ and $\ell$ is the class of a line, and $g=0$, $\overline{M}_{0,3d-1}(\mathbb{CP}^2, d\cdot\ell)$ can be thought of as the moduli stack parameterizing nodal, rational plane curves of degree d with $3d-1$ marked points. 

Diverging slightly from the ethos of previous sections, in order to motivate certain quadratically enriched results we use the remainder of this section to define Gromov--Witten invariants. Defining Gromov-Witten invariants is best done in cases, and we focus on the case when the relevant moduli space is smooth. When $X$ is convex in the sense of \cite[2.4.2]{kontsevich-manin}, for example $X=\mathbb{CP}^r$ for some $r$, $\MzeronbarXb$ is smooth \cite[1.3.2]{kontsevich1}. For each $1\leq i\leq n$ there is an evaluation map 
\[
\ev_i\colon \MzeronbarXb \to X, \hspace{5mm} (C\stackrel{f}{\to}X, p_1, \ldots, p_n) \mapsto f(p_i),
\]
and these can be defined more generally when $X$ is not convex and $g>0$. Given $p_1, \ldots, p_n$ be $n$ points of $X$ and their Poincar\'{e} duals by $\gamma_i\in H^{2\text{dim}(X)}(X)$, the cohomology class $\ev_i^*\gamma_i\in H^*(\MzeronbarXb)$ represents the class of stable maps whose image passes through $p_i$ for each $i$. Thus the $n$-fold cup product $\ev_1^*\gamma_1 \cup \cdots \cup \ev_n^*\gamma_n$ in $H^*(\MzeronbarXb)$ represents the class of stable maps whose image passes through $p_i$ for \emph{all} $1\leq i\leq n$.  The \emph{Gromov--Witten invariant counting genus $0$, degree $\beta$ stable maps passing through $p_1, \ldots, p_n$} is the degree
\begin{equation}\label{eq:smooth_GW}
\langle p_1, \ldots, p_n\rangle_{g, \beta}^X := \text{deg}\left([\MzeronbarXb] \cap (\cup_{i=1}^n \ev_i^*\gamma_i)\right).
\end{equation}
There are more general ways of defining Gromov--Witten invariants, but this formulation of the definition is most relevant for \Cref{section:quadratically_enriched}. For instance, we can define Gromov--Witten invariants counting stable maps passing through cycles $V_i$ with Poincar\'e duals $\gamma_i\in H^{n_i}(X)$ using the degree formulation in \eqref{eq:smooth_GW} more generally so long as $[\MzeronbarXb] \cap \cup_{i=1}^n \ev_i^*\gamma_i$ is a zero cycle. 

\begin{remark}\label{rmk:VFC} An immediate consideration in the definition above arises when $X$ is not convex or when $g>0$, in which case $\MgnbarXb$ is not smooth. In these instances, we cannot hope to take the degree of $\ev_1^*\gamma_1 \cup \cdots \cup \ev_n^*\gamma_n$ by capping with the fundamental class $[\MgnbarXb]$ and pushing forward to a point, as $\MgnbarXb$ is singular and has irreducible components of varying dimensions. Beautiful work of Li-Tian \cite{LT-VFC} and Behrend-Fantechi \cite{BF-VFC} shortly thereafter constructs a \emph{virtual fundamental class} for finite type Deligne-Mumford stacks. In particular, there is a virtual fundamental class $[\MgnbarXb]^{\vir}_{E^\bullet}$ of the expected dimension given a perfect obstruction theory $E^\bullet$ for $\MgnbarXb$ \cite{behrend-GW}. The virtual fundamental class $[\MgnbarXb]^{\vir}_{E^\bullet}$ can be used to define Gromov--Witten invariants in general. Gromov-Witten invariants have been used to answer Problem \ref{problem:rational_curve_count} in various cases, two of which are discussed below. 
\end{remark} 

Based on the explanations given thus far, which are far from thorough, the definition of Gromov-Witten invariants  leaves much to the imagination in terms of computational feasibility. Stunningly, powerful tools make these invariants computable in a number of cases of interest. We give two examples that will be relevant in \Cref{section:quadratically_enriched}.

\noindent {\bf Gromov--Witten invariants of blow-ups of projective space. } A natural course of study following the curve counting results of Caporaso-Harris \cite{caporaso-harris} and Kontsevich-Manin \cite{kontsevich-manin} on enumerative curve counts for $\mathbb{P}^n$ is the study of Gromov--Witten invariants of smooth, projective rational surfaces, i.e., surfaces which are deformation equivalent to $\mathbb{P}^1\times \mathbb{P}^1$ or a blow-up of $\mathbb{P}^2$ at finitely many points $x_1, \ldots, x_r$. Gromov--Witten invariants of blow-ups of projective space have been studied by G\"{o}ttsche-Pandharipande and Gathmann \cite{gottsche-pandharipande, gathmann} amongst others; Gathmann studies the case of blow-ups of $\mathbb{P}^n$. In particular, such Gromov--Witten invariants can be interpreted in terms of counts of rational curves in $\mathbb{P}^2$ with specified tangent multiplicities at the points $x_1, \ldots, x_r$. 

\noindent {\bf Equivariant localization. }  Equivariant localization is a powerful tool in differential and algebraic topology 
which allows one to compute integrals over a suitable space with a given group action in terms of the fixed point locus of the group action. Building on results of Atiyah-Bott in equivariant cohomology \cite{AB-localization} and Edidin-Graham in equivariant Chow groups \cite{EG-localization}, Graber-Pandharipande give a $\mathbb{C}^*$-equivariant virtual localization formula for a $\mathbb{C}^*$-equivariant scheme $X$ with a $\mathbb{C}^*$-equivariant perfect obstruction theory $E^{\bullet}$, which expresses $[X]^{\vir}_{E^\bullet}$ in terms of the virtual fundamental classes and Euler numbers of the $\mathbb{C}^*$-fixed point locus \cite{GP-virtual-localization}. Graber-Pandharipande show a consequence of this is a virtual localization formula for $\MgnbarXb$, which can be leveraged to express Gromov--Witten invariants $X = \mathbb{CP}^r$ as a sum over graphs corresponding to the $\mathbb{C}^*$-fixed point loci \cite{GP-virtual-localization}.

These example cases motivate the work to-date on counting rational curves in quadratically enriched enumerative geometry, covered in \Cref{section:quadratically_enriched}. To conclude this section, we give a brief overview of the specific historical developments since the advent of Gromov-Witten theory that most influence current work in quadratically enriched enumerative geometry today.  

While powerful and quite general, Gromov--Witten invariants are not perfect. They are often not enumerative, for example for threefolds in $g>0$ cases they are typically rational numbers rather than integers to account for automorphisms of stable maps. One of the most prominent paths toward resolving some of the difficulties of Gromov--Witten invariants is the introduction of Donaldson--Thomas invariants, which seek to count stable sheaves in a given curve class on a Calabi-Yau 3-fold \cite{DT1, thomas1}. The moduli space of such sheaves has a perfect obstruction theory in the sense of Behrend-Fantechi \cite{BF-VFC} that is in fact symmetric \cite{behrend-DT, BF-DT}. Donaldson--Thomas invariants can be defined by integrating over the associated virtual fundamental class of the moduli of stable sheaves. While Donaldson--Thomas invariants have their own obstacles, see  \cite{13/2-notes} for a discussion, they are integral and can be computed motivically using the Behrend constructible function on the moduli space \cite{behrend-DT}. The comparison between Gromov--Witten invariants and Donaldson--Thomas invariants is natural to explore, with several groundbreaking results on their connections \cite{MNOP1, MNOP2, MNOP-BP, MNOP-OP, MNOP-pardon}. See Thomas-Pandaripande \cite{13/2-notes} for a description of curve counting theories more generally, beautifully elucidating how Gromov--Witten and Donaldson--Thomas invariants fit into broader enumerative theories. 

In addition to counting stable maps to a target variety, significant study has also been devoted to counting curves in a fixed linear system. Notably is the Yau--Zaslow formula of Bryan--Leung, which counts genus $g$ curves with $n$ nodes in a fixed linear system on a K3 surface, remarkably given by the Dedekind $\eta$ function \cite{yau-zaslow}. See also work of Beauville, Chen, and Fantechi-G\"{o}ttsche-van Straten \cite{beauville-yz, chen-K3, fantechi-gottsche-vanstraten}. We also mention the G\"{o}ttsche conjecture, first proved by Kool-Shende-Thomas and separately Tzeng shortly after \cite{KST, tzeng}, counting nodal curves in a linear series defined by a sufficiently ample line bundle $L$ on a smooth projective complex surface $S$ in terms of a polynomial in $L^2$, $L.K_S$, $K_S^2$, and $c_2(S)$. This line of questioning around counts of nodal curves in a fixed linear series has been studied in quadratically enriched enumerative geometry, results will be mentioned in Section \ref{section:quadratically_enriched}.

As \Cref{subsection:quadratic_GW} focuses on quadratic enrichments of Gromov--Witten invariants, we end this section by noting that enriched Gromov--Witten invariants have appeared in existing work. For example, celebrated work of Y-P Lee constructs Gromov--Witten invariants in $K$-theory via integration over a virtual structure sheaf \cite{YP-Lee-GW}, and Gu\'{e}r\'{e} recently studied $K$-theoretic Gromov--Witten invariants using virtual localization for a finite group action \cite{guere-KGW}.

\section{Quadratically enriched enumerative geometry}\label{section:quadratically_enriched}

At the beginning of the 21st century, \emph{motivic homotopy theory} emerged as a popular exciting new direction in mathematics, rising to the forefront after Voevodsky's Fields Medal-winning resolution of the Bloch-Kato conjectures. Following work of Morel, it was understood that motivic spaces, the main objects of study in motivic homotopy theory, over a field $k$ admit a \emph{quadratic Euler characteristic} valued in the Grothendieck--Witt ring $\GW(k)$. Following Marc Hoyois' thesis \cite{Hoyois}, which explored further applications of this $\mathbb{A}^1$-enhancement of the Euler characteristic, Marc Levine, and simultaneously Jesse Kass and Kirsten Wickelgren, began exploring potential ways to build an enumerative geometry whose answers take values in $\GW(k)$ rather than $\Z$. This is now called $\mathbb{A}^1$\emph{-enumerative geometry} or sometimes \emph{quadratically enriched enumerative geometry}.

\begin{definition}\label{def:GW}
The \textit{Grothendieck--Witt ring} $\GW(k)$ of a field $k$ is the group completion of the semi-ring of isomorphism classes of non-degenerate symmetric bilinear forms over $k$. Explicitly it is generated by the rank one forms
\begin{align*}
    \left\langle a \right\rangle \colon k \times k &\to k \\
    (x,y) &\mapsto xay.
\end{align*}
\end{definition}
In 2000, Barge and Morel developed a theory of \emph{oriented Chow groups} (also called \emph{Chow--Witt groups}) , which are twisted by line bundles over the input scheme, and which have a natural home in the world of motivic homotopy theory. These should be thought of as an enhancement of Chow groups, in that they are generated by cycles equipped with some extra ``algebraic orientation data'' which we neglect to define precisely here. The properties of these groups were established in further detail by Fasel and Srinivas \cite{fasel_groupes_2008,fasel_chowwitt_2009}, and it is natural to think of them as providing an enhanced setting for intersection theory, where the degree takes values in the Grothendieck--Witt ring. This is the primary tool leveraged by Levine to explore quadratically enriched enumerative geometry \cite{levine_aspects_2020}.

Analogous to how the top Chern class provides an enumerative count of the zeros of a section of a corank zero vector bundle along a smooth compact manifold (or how the \emph{absolute Euler class} does the same for relatively oriented bundles in the real setting), one has an \emph{Euler class} for relatively oriented vector bundles in the motivic setting. More precisely:

\begin{definition}\label{def:euler-class-GW}
Let $X$ be a smooth proper $k$-scheme of dimension $n$, and let $V\to X$ be an algebraic vector bundle of rank $n$, which is \emph{relatively oriented}, in the sense that there is an isomorphism of line bundles
\begin{equation}\label{eqn:label}
\begin{aligned}
    \Hom(\det TX, \det V) \cong \mathcal{L}^{\otimes2}
\end{aligned}
\end{equation}
for some line bundle $\mathcal{L} \to X$. Then this bundle admits a well-defined \emph{Euler number} $n(V) \in \GW(k)$.
\end{definition}

Roughly speaking, a relative orientation can be thought of as a way of choosing a square root of the Jacobian determinant bundle, allowing one to quadratically count solutions to enumerative problems arising from Euler number computations in a way that is analogous to keeping track of signs in real enumerative geometry (see subsection \ref{section:real-signed-counts}. The remarkable result is that the Euler number provides a quadratically enriched count of the zeros of an algebraic section of the bundle in a way that is independent of the choice of section. This result can be found in \cite[1.1]{bachmann_euler_2023}, but is the culmination of a lot of work towards establishing a rigorous theory of motivic Euler classes \cite{BargeMorel,Morel,asok_comparing_2016,levine_aspects_2020,LevineRaksit,KW27}.

A seminal application of these techniques is the following quadratically enriched count of 27 lines on a cubic surface.

\begin{theorem}[{\cite{KW27}}]
\label{thm:enriched-27}
We have that
\begin{align*}
    n \left( \Sym^3\mathcal{S}^\ast \to \Gr_k(2,4) \right) = 15 \left\langle 1 \right\rangle + 12 \left\langle -1 \right\rangle\in\GW(k).
\end{align*}
\end{theorem}
The rank of this form is 27, recovering the classical count of 27 lines on a smooth cubic surface. The signature is 3, which recovers Segre's theorem discussed in \Cref{section:real}. Over finite fields, for examples, this reveals new constraints on the possible fields of definition for hyperbolic and elliptic lines on cubic surfaces. This result has been generalized to provide a count of lines on symmetric hypersurfaces in general (see \cite{Pauli-quintic3fold}, \cite[\S8]{levine_motivic_2019}, and \cite[\S6.1]{bachmann_euler_2023}).

\subsection{Computing local indices} One of the core problems in quadratically enriched enumerative geometry is computing local indices. As the Euler class is the sum of the local indices of the zeros of a section, the local index should be read as some local geometric data that a single solution contributes to the overall count (think: a line being hyperbolic versus elliptic).

In the setting of \Cref{def:euler-class-GW}, the local index is computed as an $\mathbb{A}^1$-\emph{Brouwer degree} of the intrinsic derivative of the section $\sigma \colon X \to V$ around an isolated zero. One can pass to affine charts (see \cite[Lemma~19]{KW27} for precise details), hence the problem reduces to producing an element in $\GW(k)$ from an endomorphism of affine $n$-space with an isolated zero at the origin. In \cite{KW-EKL}, the authors argued that the $\mathbb{A}^1$-Brouwer degree at a rational point generalizes the \emph{Eisenbud--Khimshiashvili--Levine signature formula} over the reals \cite{eisenbud_algebraic_1977,himsiasvili_local_1977}. This form can be interpreted in a number of ways -- it is a quadratic duality pairing explored by Scheja and Storch for complete intersections \cite{scheja_uber_1975}, or can be thought of more generally as a trace form arising from coherent duality. The so-called \emph{EKL form} is very computable, and was shortly generalized to local degrees at points with residue fields finite separable over the base \cite{brazelton_trace_2021} and finally points with arbitrary residue field \cite{brazelton_bezoutians_2023}. It can now be computed explicitly in \texttt{Macaulay2} \cite{borisov_bbb_2024}.

\subsection{More results in quadratically enriched enumerative geometry} Since the seminal work of Kass-Wickelgren and Levine, a tremendous amount of progress has been made leveraging these techniques to provide enriched counts of classical questions in enumerative geometry. We highlight a few of these which are not mentioned elsewhere in this paper.

 A quadratically enriched Schubert calculus was developed by studying the Chow--Witt groups of Grassmannians and flag varieties \cite{wendt_chow-witt_2024,hudson_chow-witt_2024}. Many of these results are obtained as an $\mathbb{A}^1$ degree or Euler number. B\'ezout's theorem has seen a quadratic enrichment \cite{mckean_arithmetic_2021}, as has the Circles of Apollonius problem explored in \Cref{problem:apollonius} \cite{mckean_circles_2022}. An enriched count of bitangents to a planar curve has been developed \cite{LarsonVogt} and was further explored in \cite{kummer_bounding_2024}. An enriched count of lines meeting four lines in three-space \cite{srinivasan_arithmetic_2019} and higher-dimensional analogues of this \cite{brazelton_enriched_2025} have been developed. Arithmetic inflection for linear series along curves has been enriched \cite{cotterill_arithmetic_2023,cotterill_arithmetic_2024}. 
 Other results include \cite{agostini_ulrich_2025,darwin_conics_2023,darwin_quadratically_2022,espreafico_motivic_2023,espreafico_quadratic_2025,kim_global_2024, muratore_arithmetic_2025}. Kummer uses quadratically enriched invariants to give a signed count of 2-torsion points on real abelian varieties over $\mathbb{R}$ and $\mathbb{C}$ \cite{kummer_signed_2023}. 
 Pajwani-P\'{a}l make meaningful progress toward an arithmetically enriched Yau-Zaslow formula in characteristic 0 in \cite{pajwani_arithmetic_2025}, which is an example of a quadratically enriched rational curve count in a linear series. While this is different from the quadratically enriched rational curve counts through given point conditions, this naturally leads to the next subsection. 

\subsection{Quadratically enriched curve counting} \label{subsection:quadratic_GW}

A problem that has appeared throughout this article culminates in the following question:  

\begin{problem}
    Let $k$ be a field. Given general points $p_1, \ldots, p_r$ of $\mathbb{P}^2_k$ such that all residue fields $k(p_i)$  are separable over $k$, how many degree $d$ rational plane curves pass through $p_1, \ldots, p_r$? 
\end{problem}

A quadratically enriched answer to this question would give a quadratic form whose rank is equal to the Gromov--Witten invariant $N_d$, recovering the complex count of rational degree $d$ plane curves when $r=3d-1$, and whose signature is equal to the Welschinger invariant $W_{d}$, recovering the signed count of real degree $d$ plane curves. Recent work of J. Kass, M. Levine, J. Solomon, and K. Wickelgren introduce precisely these quadratically enriched counts of rational curves in a certain divisor class on a del Pezzo surfaces of degree $\geq 4$ over perfect fields of characteristic $\neq 2, 3$ \cite{KLSW-rational-curves}. This builds on prior work of M. Levine defining quadratically enriched Welschinger invariants \cite{levine-welschinger}. Beyond recovering real and complex rational curve counts, their work proposes a definition of rational curve counts over other fields. 

Let $X$ denote a del Pezzo surface of degree $\geq 4$ over a perfect field of characteristic $\neq 2,3$. Using the notation of \Cref{section:Gromov-Witten}, let $\overline{M}_{0,n}(X,d)$ denote the moduli space of genus 0 stable maps of degree $d$ to $X$, and let 
\[
\ev\colon \overline{M}_{0,n}(X,d) \to X^n, \hspace{5mm} (C\stackrel{f}{\to}X, p_1, \ldots, p_n) \mapsto (f(p_1),\ldots, f(p_n))
\]
denote the total evaluation map. Kass-Levine-Solomon-Wickelgren \cite{KLSW-rational-curves} define the quadratically enriched rational curve counts to be 
\[
\text{N}_{X, d, \sigma}:= \deg^{\mathbb{A}^1}(\ev),
\]
analogously to the definition in \Cref{eq:smooth_GW} defining Gromov--Witten invariants as a degree. Examples are given in \cite[Table 1]{KLSW-rational-curves}. This definition relies on the technical condition of relative orientability of the total evaluation map, which is shown to be relatively oriented away from a high codimension locus in $X$ in \cite{KLSW-orientation}. In general, $\text{deg}^{\mathbb{A}^1}(ev)$ is a quadratic form over $X^n$, not over $k$, but a unique class in $\GW(k)$ can be obtained when $X$ is $\mathbb{A}^1$-connected, for example when $X$ is a smooth proper rational surface. 

Remarkably, the quadratic forms $\text{N}_{X,d,\sigma}$ depend only on the list of separable field extensions over $k$ determined by the chosen points, $\{k(p_1), \ldots, k(p_r)\}$ \cite[Section 8]{KLSW-rational-curves}. Equally remarkably, $\text{N}_{X,d,\sigma}$ can be computed as a sum of individual contributions of each curve, which are quadratic enrichments of the Welschinger mass (\Cref{def:welschinger-invariant}), which take into account the fields of definition of branches of the curve at nodal points. In some cases, the numbers $\text{N}_{X,d,\sigma}$ can be computed using the $\mathbb{A}^1$-Euler number of a certain oriented vector bundle over $X$, see \cite[9.1]{KLSW-rational-curves}. 
 
This work has generated significant activity since its inception. This includes connections with tropically enriched enumerative geometry, which is discussed in the next subsection. See also \cite{brugalle-wickelgren, chen-wickelgren}. 

Importantly, this work occurs alongside the study of quadratically enriched Donaldson--Thomas invariants. Given a smooth, projective 3-fold $X$, Viergever and Viergever-Levine define quadratically enriched Donaldson--Thomas invariants using a perfect obstruction theory for the Hilbert scheme $\text{Hilb}^n(X)$ of ideal sheaves of length $n$ with support of dimension 0 on $X$ \cite{quadratic-DT2, quadratic-DT1}. There is an associated virtual fundamental class for $\text{Hilb}^n(X)$ using the construction of \cite{levine-VFC}, which is used to define the enriched Donaldson--Thomas invariants. Viergever and Viergever-Levine leverage work of Levine on virtual localization and a relative orientation for $\text{Hilb}^n(X)$ for computations, see \cite{levine-localization-witt, levine-virtual-localization,levine-orient-hilbert}. Work of Espreafico-Walcher in \cite{quadratic-DT3} leverages the realization of $\text{Hilb}^n(\mathbb{A}^3)$ as the critical locus of a function and the compactly supported $\mathbb{A}^1$-Euler characteristic to define quadratically enriched Donaldson--Thomas invariants for $\mathbb{A}^3$ in $\GW(k)$. 

\subsection{Tropical and quadratic enumerative geometry} A well-established direction in enumerative geometry (which we do not explore in these notes) is the use of methods from \emph{tropical geometry} to carry out computations in classical enumerative geometry. A beautiful example of this is the \emph{Mikhalkin correspondence}, which states that the number $N_d$ of rational degree $d$ curves through $3d-1$ generic points on $\mathbb{P}^2$ can be computed as a count of rational tropical curves with a given Newton polytope through $3d-1$ generically chosen points on the plane (computed with multiplicity) \cite{mikhalkin_enumerative_2005}. Mikhalkin further proved an analogue for real rational curves interpolating points, which provides a tropical way to compute the signed count of these curves weighted by their Welschinger invariants.

With this and $\mathbb{A}^1$-enriched Welschinger invariants in mind, it is natural to ask whether these tropical methods developed by Mikhalkin can be given a quadratic enrichment in order to capture both the classical and real counts simultaneously. This leads to a beautiful sub-field of quadratically enriched enumerative geometry which can leverage tropical methods to compute enriched solutions \cite{markwig_bitangents_2023,jaramillo_puentes_arithmetic_2024,jaramillo_puentes_quadratically_2025}. For an expository introduction to these ideas see \cite{pauli_pcmi_2024}.

\section{A few emergent directions}

There are many exciting new directions in enumerative geometry, and far too many to do justice in this paper. Here we highlight two directions, one coming from the theory of \textit{random algebraic geometry} and the other coming from equivariant mathematics.

\subsection{Random enumerative geometry}

Over the reals, conservation of number fails as previously discussed. However we can still ask the question: what is the \emph{expected} number of real solutions? Here we discuss two ways to approach such a problem, one coming from probability theory and the other from Hodge theory and hyperbolic geometry.

The easiest place to see the failure of conservation of number (and hence the easiest place to establish a testbed for a random theory of enumerative geometry) is in the fundamental theorem of algebra --- that a univariate polynomial $f\in \R[x]$ of degree $n$ may fail to have $n$ real roots, counted with multiplicity. Motivated by this question, Kac investigated the following question in the 1940's: what is the \emph{expected number} of real roots of a randomly chosen real polynomial of some degree $n$? One first has to define what they mean by ``random,'' and in Kac's work he assumes that the coefficients are distributed according to a normal distribution on $\R$. With this in mind he gives an exact formula for the expected value, together with an asymptotic of $\frac{2}{\pi}\log(n)$ \cite{kac_average_1942,kac_average_1948}. An entirely different, although spiritually similar, approach to this problem came from the work of Stephen Oswald Rice, a researcher at Bell Labs working in signal processing and random noise \cite{rice_mathematical_1944}. These results can be extended to provide formulas for random maps on manifolds \cite[\S4]{breiding_random_nodate} -- these sorts of generalizations are called \emph{Kac-Rice formulas}.

This flavor of question (what is the expected behavior of a randomly chosen algebraic object over the reals) leads to a new program of mathematics called \emph{random algebraic geometry} \cite{lerario_what_nodate}. A key application of these techniques is to leverage tools like the Kac-Rice formula to compute the expected number of real solutions to an enumerative problem -- we may call this \emph{random enumerative geometry}. A motivating result in this direction is the following:

\begin{theorem}[{\cite{basu_random_2019}}]
The expected number of real lines on a random real cubic surface is $6\sqrt{2}-3$.
\end{theorem}

This theorem is proven by choosing random sections of the symmetric bundle $\Sym^3\mathcal{S}^\ast$ over the Grassmannian of lines in $\RP^3$, then running the aforementioned machinery to compute the expected number of zeros. Implicit in this work is a choice of a certain Gaussian probability distribution, and the answer to the problem will very naturally depend on this choice. These methods have been used, among other applications, to study random Schubert calculus \cite{burgisser_probabilistic_2020} and to extend these randomized questions away from the reals and toward the $p$-adics \cite{ait_el_manssour_probabilistic_2022}.

Returning to lines on real cubic surfaces, an entirely different approach to the same problem above comes from analyzing the moduli of real cubic surfaces, and attempting to compare the volumes of the five connected components. A priori this is a poorly phrased problem, as there is no natural metric to compute volume on the moduli space of cubic surfaces (being, as it is, a quotient of a Zariski open subset of projective space by the action of an infinite group). Miraculously, we can leverage techniques from Hodge theory to endow it with a hyperbolic metric! More precisely, we can study the variation of Hodge structure on real cubic threefolds covering our real cubic surfaces, analyze their images under a period map, and finally compute the volume of each connected component of the moduli space via Vinberg's algorithm or similar techniques. For the specific problem of computing lines on real cubic surfaces, this was accomplished in groundbreaking work of Allcock-Carlson-Toledo \cite{allcock_hyperbolic_2010}. This program of mathematics has already seen wide-ranging applications, from the classification of real cubic fourfolds \cite{finashin_topology_2010} to the study of real binary octics \cite{chu_geometry_2011}.

As to the potential of both approaches for gaining intuition and revealing structure in real enumerative geometry, it seems we have only just scratched the surface of what is possible.

\subsection{Equivariant enumerative geometry}\label{subsec:random}
In this last subsection, we abandon any pretense of modesty to discuss an ongoing program of work by the authors and others to leverage tools from equivariant homotopy theory in order to provide equivariantly enriched counts of classical questions in enumerative geometry. Given an enumerative problem and some symmetry, it is natural to ask how the symmetry group interacts with the solutions to the problem.\footnote{Related but orthogonal work in this direction includes \cite{roberts_equivariant_1985,damon_g-signature_1991,costenoble_algebraic_2024}} 

When $G$ is a finite group, we can build a theory of equivariant Euler classes for $G$-equivariant complex topological vector bundles along $G$-manifolds, which allows us to prove an \emph{equivariant conservation of number} result under mild hypotheses \cite{EEG}. Roughly speaking this states that the symmetries of an equivariant enumerative problem are always conserved. For instance, given a smooth cubic with automorphism group $G$, the group always acts on the 27 lines in the same way \cite{EEG}. Another application of the equivariant Euler number has been given in \cite{BB-bitangents}, where the authors give an enriched count of orbits to the 28 bitangents of any smooth, non-hyperelliptic, symmetric quartic curve. Equivariant counts of orbits of solutions to enumerative problems have also appeared in the context of counting solutions to enumerative problems in families that are invariant under a finite group action. Equivariant counts have been provided for counting nodes in a $G$-invariant pencil of conics \cite{bethea_enriched_2025} and rational cubics interpolating a $G$-invariant set $8$ general points in $\mathbb{CP}^2$ \cite{bethea_equivariant_2025}, both of these results recover a real signed count of nodal conics and rational cubics respectively when $C_2$ acts on $\mathbb{CP}^2$ by pointwise complex conjugation. The intersection of real symmetric hypersurfaces has also been explored \cite{lidz_intersections_2024}. Equivariantly enriched Gromov-Witten invariants for smooth, projective complex varieties with the action of a finite group will appear in upcoming work of the first named author and Wickelgren. 

An interesting direction is to ask how symmetry interplays with the Galois group of an enumerative problem (in the sense of \cite{Hermite51,Jordan,Harris-Galois}). Using variation of Hodge structure techniques analogous to those discussed in \Cref{subsec:random}, the second-named author and Raman showed the Galois group of lines on a $S_4$-symmetric cubic surface is the Klein four-group \cite{brazelton_monodromy_2025}. Using entirely different techniques derived from the world of stacks, Landi investigated the same question for cubic surfaces with an involution and computed the Galois group of their lines \cite{landi_stacks_2025}, and Pichon-Pharabod and Telen extended this work by numerically certifying monodromy computations of Galois groups for all symmetric cubic surfaces \cite{pichonpharabod2025galoisgroupssymmetriccubic}. Upcoming work of the second named author, Landi and Raman explores these ideas in greater detail.

\bibliographystyle{amsalpha}
\bibliography{citations}
\end{document}